\documentclass[10pt]{amsart}

\usepackage{amsmath,amssymb,amsfonts, amscd,verbatim}
\usepackage[mathscr]{eucal}
\usepackage{diagrams}          

\title[Correction to ``Duality and Flat Base Change \dots'']
{Correction to the paper\\ 
   ``Duality and Flat Base Change on\\
           Formal Schemes''}

\author[L.\;Alonso]{\smash{Leovigildo Alonso Tarr\'{\i}o}}
\address[L.\;A.\;T.]{Departamento de \'Alxebra\\
Facultade de Matem\'a\-ticas\\
Universidade de Santiago de Compostela\\
E-15782  Santiago de Compostela, SPAIN}
\email{leoalonso@usc.es}

\author[\smash{A.\;Jerem\'{\i}as}]{\smash{Ana Jerem\'{\i}as L\'opez}}
\address[A.\;J.\;L.]{Departamento de \'Alxebra\\
Facultade de Matem\'a\-ticas\\
Universidade de Santiago de Compostela\\
E-15782  Santiago de Compostela, SPAIN}
\email{jeremias@usc.es}

\author[J.\;Lipman]{Joseph Lipman}
\address[J.\;L.]{Dept.\ of Mathematics\\
  Purdue University\\
  W. Lafayette IN 47907, USA}
\email {lipman@math.purdue.edu}

\thanks{First two authors partially supported by Spain's DGESIC PB97-0530
research project. They thank the Mathematics Department of Purdue
University for its hospitality and support.}

\thanks{Third author partially
             supported by the National Security Agency.\vadjust{\kern1.5pt}}

\theoremstyle{plain}
\newtheorem{thm}{Theorem}
\newtheorem{lem}[thm]{Lemma}

\newtheorem{prop}[thm]{Proposition}

\theoremstyle{remark}

\theoremstyle{definition}

\DeclareMathOperator{\Hom}{Hom}
\DeclareMathOperator{\spf}{Spf}

\newcommand{\E}{{\mathcal E}}
\newcommand{\F}{{\mathcal F}}
\newcommand{\G}{{\mathcal G}}
\newcommand{\cH}{\mathcal H}
\newcommand{\I}{{\mathcal I}}

\newcommand{\CO}{{\mathcal O}}

\newcommand{\D}{{\mathbf D}}
\newcommand{\PP}{{\mathbf P}}
\newcommand{\R}{{\mathbf R}}

\newcommand{\sHomb}{\cH{om}^{\bullet}}
\newcommand{\iGp}[1]{{\varGamma_{\<\!#1}'}}
\newcommand{\ush}[1]{{#1^{\textup{\texttt\#}}}}

\newcommand{\Rfs}{{\mathbf R f_{\!*}}}
\newcommand{\comp}{\boldsymbol\Lambda}
\newcommand{\smcirc}   
  {{\raise.15ex\hbox to.7em{$\hss \scriptstyle\circ\hss$}}} 
\newcommand{\lc}{\cite{AJL}}
\newcommand{\<}{\mspace{-1mu}}
\renewcommand{\>}{\mspace {1mu}}
\newcommand{\st}{\scriptstyle }

\newcommand{\X}{\mathscr{X}}
\newcommand{\Y}{\mathscr{Y}}
\newcommand{\U}{\mathscr{U}}
\newcommand{\V}{\mathscr{V}}
\renewcommand{\P}{\mathscr{P}}

\newcommand{\Dc}{\D_{\mkern-1.5mu\mathrm c}}
\newcommand{\Dqc}{\D_{\mkern-1.5mu\mathrm {qc}}}
\newcommand{\wDqc}{ \widetilde
         {\vbox to6.5pt{\vss\hbox{$\mathbf D$}}}      
          _{\mkern-1.5mu\mathrm {qc}} }
\newcommand{\wDqcp}{\wDqc^{\lower.5ex\hbox{$\scriptstyle+$}}}
\newcommand{\Dqcp}{\Dqc^{\lower.5ex\hbox{$\scriptstyle+$}}}

\newcommand{\lto}{\longrightarrow}

\newcommand{\liso}{\overset{\sim}{\lto}}

\newcommand{\iso}%
{{\mkern8mu\longrightarrow \mkern-25.5mu{}^\sim\mkern17mu}}
\newcommand{\osi}%
{{\mkern8mu\longleftarrow \mkern-24.5mu{}^\sim\mkern16mu}}
\newcommand{\Iso}{\vbox to
        0pt{\vss\hbox{$\widetilde{\phantom{ne}}$}\vskip-7pt}}

\newcommand{\under}[2]
{\vbox to 0pt{\vskip-#1 ex\hbox{$\scriptstyle #2$}\vss}}
\newcommand{\dirlm}[1]%
  {
     {\lim\hskip-1.58em\lower.65ex
       \hbox{$
                {}_{\stackrel{\lower1ex\hbox
                                       {$\scriptstyle -\!\!\!\longrightarrow$}    
                                      }{\vbox to0pt{\vss\vskip.6ex
                                            \hbox{$\scriptstyle{}^{#1}$}\vss}}
                   }
            $}
     }
\:}

\def\bilap#1{\hbox to 0pt{\hss#1\hss}}
 \def\Rarrow#1{\bilap{\hbox to#1{\rightarrowfill}}}

\begin{document}

\begin{abstract}
In \S8.3 of our paper ``Duality and Flat Base Change on Formal Schemes" \cite{AJL}
some important results concerning localization of, and preservation
of coherence by, basic duality functors,
were based on the false statement that any closed
formal subscheme of an open subscheme of the completion $\mathscr P$ of a
relative projective space is an open subscheme of a closed formal subscheme of
$\mathscr P$. In this note, the said results are provided with solid foundations.
\end{abstract}

\renewcommand{\subjclassname}{Mathematics Subject Classification
\textup{2000}}

\maketitle

In Proposition 8.3.1 of our paper \cite{AJL},  the duality functors 
$f^!$ and $\ush f$
associated to a pseudo-proper map~$f\colon \X\to \Y$  of noetherian formal
schemes (i.e., right adjoints of  suitable restrictions of the derived 
direct-image
functor $\Rfs$) are asserted to be local on
$\X$, as a consequence of flat base change. Moreover,  in Proposition 
8.3.2 it is
asserted that (roughly speaking)
$\ush f$ preserves coherence. Brian Conrad pointed out that our 
justifications are
deficient because they use the claim 8.3.1$\>$(c) that a map between
noetherian formal
schemes that can be factored as a closed immersion followed by an open one can
also be factored as an open immersion followed by a closed one, which 
is not true in general.%
\footnote{The ``proof" breaks down in the second-last line of
\cite[p.\,88]{AJL}, where it is erroneously stated that localization 
followed by
completion commutes with forming kernels of  homomorphisms of adic
rings.\looseness=-1}
Indeed,  Conrad observed that for any 
$(A, x, p)$ with $A$ an adic domain, $x\in A$ such that $B\!:=A_{\{x\}}$ is a domain,
and 
$p$  a nonzero $B$-ideal contracting to (0) in $A$, the natural map 
${\rm Spf}(B/p)\to{\rm Spf}(A)$ is a counterexample. Such a triple was
provided to us by Bill Heinzer:\vspace{1.5pt}

\begin{footnotesize}
With $w,x,y,z$  indeterminates over a field $k$, set
$$
A\!: = k[w,x,z][[y]] \quad\textup{and}\quad 
B\!:  = A_{\{x\}}=k[w, x,1/x,z][[y]].
$$  
Let $P$ be the prime ideal $(w,z)A$ and 
$R\!:=A_P\subset B_{P\<B}=\colon \!S$, so that  
$R\subset S$ are 2-dimensional regular local domains such that 
the residue field
of~$S$ (i.e., the fraction~field of~$k[x, 1/x][[y]]$) is
transcendental over that
of~$R$ (i.e., the fraction~field of~$k[x][[y]]$).
Then \cite[p.\,364, Theorem 1.12]{HR} says that 
\emph{there exist  infinitely many height-one prime $S$-ideals
in the generic fiber over $R$.} 

Any of these contracts in $B$ to a (prime) $p$ as
above.\par
\end{footnotesize}

\smallskip
Our purpose
here is to validate the aforementioned Propositions
by means of a ``localization" Lemma (\ref{fix} below). Thus all  other results in
\lc\  depending on these Propositions remain as they are. (No other results depend
on the faulty 8.3.1$\>$(c).)\looseness=-1

However,
Proposition 8.3.1 is weakened in that we get  isomorphisms which are not
\emph{a priori} functorial or canonical at the level of derived 
categories but only
at the level of homology sheaves. This drawback does not affect the
applications. 

(Derived functoriality\- and canonicity might well be attainable, for
example through a suitable variant---if such exists---of compactification of
separated pseudo\-finite type maps of noetherian formal schemes.)

\smallskip
Notation and terminology are as in \lc\ (which has an index starting 
on p.\,125). For example, with $\I_\Y$ an ideal of definition of the
noetherian formal scheme~$\Y$ and
$\iGp\Y(-)$ the torsion functor
$\dirlm{n>0\,}\cH{om}(\CO_\Y/\I_\Y^n,-)$, an
$\CO_\Y$- complex
$\F$ lies in
$\wDqcp(\Y)$  if the homology $H^i(\F)$ vanishes for~$i\ll0$ and 
the derived-torsion complex $\R\iGp\Y \F$ has quasi-coherent homology.
In particular, $\Dqcp(\Y)\subset\wDqcp(\Y)$ \cite[p.\,54, 5.2.10]{AJL}.

\begin{lem}\label{fix} Let\/ 
$
\U \overset{q_k^{}}{\lto}
\V_k \overset{h_k^{}\>}{\lto} \X_k \overset{f_k^{}}{\lto} \Y
\ (k=1,2)
$
be maps of noetherian formal schemes with\/ $q_k^{}$ a closed
immersion, $h_k^{}$ an open immersion, $f_k^{}$ 
pseudo-proper, and\/ $f_1h_1q_1=f_2h_2q_2$. 
Then one can define isomorphisms
\[
q_1^!h_1^*f_1^!\F \iso q_2^!h_2^*f_2^!\F\>,\qquad\quad
\ush{q_1}h_1^*\ush{f_1}\F \iso \ush{q_2}h_2^*\ush{f_2}\F
\qquad\bigl(\F \in \wDqcp(\Y)\bigr)
\]
such that the induced  homology isomorphisms are canonical and functorial.
\end{lem}

\begin{proof} Recall that $f^!$ is the notation used for
$f_{\mathrm t}^\times$ when $f$ is
pseudo-proper, and that $\comp_\X(-)\!:=\R\sHomb(\R\iGp\X\CO_\X,-)$. The
isomorphisms of functors from $\D(\Y)$ to $\D(\U)$:
\[\comp_{\U}q_k^!h_k^*f_k^! \iso \comp_{\U}q_k^!\comp_{\<\V_k}h_k^*f_k^!
\iso\comp_{\U}q_k^!h_k^*\comp_{\X_k}f_k^! \iso \ush{q_k}h_k^*\ush{f_k}
\qquad(k \in \{1, 2\})
\]
where the second is obvious and the other two are given by
  \cite[Corollary 6.1.5]{AJL},  show that it suffices to establish the first
isomorphism in Lemma~\ref{fix}.

Let $\I_{\Y} \subset \CO_{\Y}$, $\I_{\X_1} \subset \CO_{\X_1}$ and
$\I_{\X_2} \subset \CO_{\X_2}$ be ideals of definition such that
$\I_\Y\CO_{\X_1} \subset \I_{\X_1}$ and $\I_\Y\CO_{\X_2} \subset
\I_{\X_2}$. The ideals $\I_1 := \I_{\X_1}\CO_\U$, $\I_2 :=
\I_{\X_2}\CO_\U$ and $\I := \I_1 + \I_2$ are ideals of definition of the
formal scheme $\U$, possibly different. For each $n > 0\>$ let
\( u_n \colon U_n \to \U \) be the closed immersion determined by
$\I^n\<$ (so that $U_n$ is an ordinary noetherian scheme with the same
underlying topological space as $\U$, but with structure sheaf $\CO_\U/\I^n$).
The desired isomorphism results from the existence---to be shown---of a family
of isomorphisms
\begin{equation}\label{ison}
u_{n*}u_n^!q_1^!h_1^*f_1^!\F \liso u_{n*}u_n^!q_2^!h_2^*f_2^!\F
\qquad\bigl(n>0,\;\F \in \wDqcp(\Y)\bigr),
\end{equation}
compatible with the homotopy colimit triangles given by
\cite[Lemma 5.4.1, Proposition 5.2.1(a), and Example 6.1.3(4)]{AJL},
for $k \in \{1,2\}$:
\begin{equation}\label{triangles}
\bigoplus_{n > 0}u_{n*}u_n^!q_k^!h_k^*f_k^!\F
\longrightarrow\>
\bigoplus_{n > 0}u_{n*}u_n^!q_k^!h_k^*f_k^!\F \lto
q_k^!h_k^*f_k^!\F \overset{+}{\lto}
\end{equation}
For, a basic property of triangles is that such a family of
isomorphisms extends (not necessarily uniquely!) to an isomorphism between the 
``summits'' $q_k^!h_k^*f_k^!\F$.

Though the isomorphisms \eqref{ison} will be canonical, it does not follow
that their extension to the summits is. The definition of
\eqref{ison} is based on the fact, well-known, 
though not yet published in full generality, that on
the category of separated finite-type maps of arbitrary noetherian schemes,
there is a pseudofunctor~$^{\boldsymbol !}$ taking values in $\Dqcp$, 
agreeing with $f^!$ when the map $f$ is proper and with $f^*$ when $f$ 
is an open immersion,
cf.~\cite[pp.\,303--318]{De}. The aforesaid lack of canonicity
obstructs immediate extension of the pseudofunctor in question to the
category of formal schemes, an extension whose existence would give a stronger
canonical version of Lemma~\ref{fix}. (Nevertheless such an extension might
always exist for reasons as yet unknown to us.)

However the
$i$-th homology
$H^i(q_k^!h_k^*f_k^!\F\>)$  \emph{is} canonically isomorphic
to the direct limit of  $H^i(u_{n*}u_n^!q_k^!h_k^*f_k^!\F\>)$; and so  we
will have produced  \emph{canonical functorial homology isomorphisms}
$$
H^i(q_1^!h_1^*f_1^!\F\>)\iso
H^i(q_2^!h_2^*f_2^!\F\>)\qquad(i\in\mathbb Z) .
$$

 The isomorphisms (\ref{ison}) arise from applying $u_{n*}$
to the below isomorphisms \eqref{iso2}, that we describe next.
Consider the diagram
\begin{small}
\begin{equation}\label{hex}
\begin{diagram}[htrianglewidth=2.3em,tight]
  \\
&&           &             &  U_n  &            & &&
      \\
&& & \ldTo^{\raisebox{.5ex}{$\st q_{1\>n}^{}$}}(4,3) & 
\dTo^{\raisebox{-.5ex}{$\st
u_{n}^{}\!$}} &  \rdTo^{\raisebox{.5ex}{$\st q_{2\>\>n}^{}$}}(4,3)  & &&
      \\
& & & &
      \\
V_{1\>n}&  &    &          & \U &             & &&V_{2\>\>n}
      \\
&\rdTo^{\raisebox{.5ex}{\!$\st v^{}_{1\>n}$}}(2.1,1.9)&    &
\ldTo^{\raisebox{.5ex}{$\st q^{}_1\>$}} &   &
\rdTo^{\raisebox{.5ex}{$\st q^{}_2$}}    &  &
\ldTo^{\raisebox{.5ex}{$\st v^{}_{2\>\>n}$}}&
      \\
&&\V_{\<\<1}       &      &    &    &\V_{\<2} &&
      \\
& & & &
      \\
\dTo^{\raisebox{2ex}{$\st h_{\<1\>n}\<\<$}}&\raisebox{2ex}{$\lozenge\ $}&
\dTo^{\raisebox{2ex}{$\st h_{\<1}\!$}} & &    & &
\dTo_{\raisebox{2ex}{$\st \<h_{\<2}$}}&
\raisebox{2ex}{$\ \lozenge$}&\dTo_{\raisebox{2ex}{$\st h_{\<2\>\>n}$}}
      \\
&&\X_1   &  &    &   & \X_2 &&
      \\
&\ruTo_{\raisebox{2ex}{\!$\st x^{}_{1\>n}$}}(2.1,1.9)&
& \rdTo_{\raisebox{2ex}{$\st f_{\<\<1}^{}$}}(2.1,1.9) &  &
\ldTo_{\raisebox{2ex}{$\st f_{\<2}^{}$}}(2.1,1.9) &
&\luTo_{\raisebox{2ex}{$\st x^{}_{2\>\>n}$}\!\!}(2.1,1.9)&
      \\
X_{1\>n}&&           &             &\raisebox{.3ex}{$ \Y $} & 
& &&\ X_{2\>\>n}
      \\
&\rdTo_{\raisebox{2ex}{$\st f_{\<\<1\>n}^{}$}}(4,3)
& & & & & &\ldTo_{\raisebox{2ex}{$\st f_{\<2\>\>n}^{}$}}(4,3)&
      \\
&& & & \uTo^{\raisebox{2ex}{$\st y_{\<n}$}\mspace{-1.5mu}} & & &
& & & &
      \\
  &&&                       & Y_n &             & &&
      \\
\end{diagram}
\end{equation}
\end{small}%
where $x_{1\>n}$, $x_{2\>\>n}$, $y_n$ are the closed immersions given
by the ideals
$\I^n_{\X_1}$, $\I^n_{\X_2}$ and~$\I^n_{\Y}$ respectively,
the maps $f_{1\>n}^{}$ and
$f_{2\>\>n}^{}$ are induced by $f_1^{}$ and $f_2^{}$ respectively, the
subdiagrams marked by $\lozenge$ are fiber squares, and
$q^{}_{1\>n}$ and $q^{}_{2\>\>n}$ are the closed immersions induced by
$q^{}_1$ and $q^{}_2$, respectively. The outer hexagon is then a diagram of
ordinary noetherian schemes with $f_{1\>n}^{}$ and $f_{2\>\>n}^{}$ proper maps,
$h_{1\>n}$ and $h_{2\>\>n}$ open immersions and $q^{}_{1\>n}$ and $q^{}_{2\>\>n}$
closed immersions.

Use adic flat base change
\cite[Theorem 7.4]{AJL} and pseudofunctoriality \cite[Theorem 6.1(b)]{AJL} to
obtain the natural composite isomorphism 
$$
s_{1\>n}(\F\>)\colon u_n^!q_1^!h_1^*f_1^!\F
\iso q_{1\>n}^!v_{1\>n}^!h_1^*f_1^!\F
\iso q_{1\>n}^!h_{1\>n}^*x_{1\>n}^!f_1^!\F
\iso q_{1\>n}^!h_{1\>n}^*f_{1\>n}^!y_n^!\F\<,
$$
and analogously,
\[s_{2\>\>n}(\F\>)\colon u_n^!q_2^!h_2^*f_2^!\F \iso
q_{2\>\>n}^!h_{2\>\>n}^*f_{2\>\>n}^!y_n^!\F\>. \]
Using the above pseudofunctor on ordinary schemes we write
\mbox{$h_{k\>n}^*=h_{k\>n}^!\>$;} and since
$y^!_n\F\in \Dqcp(Y_n)$ \cite[p.\,59, Theorem~6.1]{AJL},  there results a
natural isomorphism
\[
r_n(\F\>)\colon q_{1\>n}^!h_{1\>n}^*f^!_{1\>n}y_n^!\F \iso
q_{2\>\>n}^!h_{2\>\>n}^*f^!_{2\>\>n}y_n^!\F\<.
\]
(This by itself is a canonical version of Lemma~\ref{fix} for \emph{ordinary}
schemes, for which $u_n$, $v_{k\>n}$, $x_{k\>n}$ and $y_n$ in
diagram~\eqref{hex} are identity maps.) 

We have then the natural functorial isomorphisms
\begin{equation}\label{iso2}
s_{2\>\>n}(\F\>)^{-1}r_n(\F\>)s_{1\>n}(\F\>)\colon
u_n^!q_1^!h_1^*f_1^!\F\iso u_n^!q_2^!h_2^*f_2^!\F
\qquad(n>0).
\end{equation}

\smallskip
Still to show is that the isomorphisms \eqref{ison} are compatible with
the triangles~\eqref{triangles}. Let $u_{n+1}^n\colon U_n\to U_{n+1}$ be the
natural closed immersion, and let 
$$
t_n\colon u_{n*}u_n^!\cong u_{n+1*}u_{n+1*}^n (u_{n+1}^n)^!u_{n+1}^!
\to u_{n+1*}u_{n+1}^!
$$ 
be the natural map.
By the definitions involved, the desired compatibility amounts to commutativity of
the following diagram in the category of functors from
$\wDqcp(\Y)$ to $\D(\U)$, where $c_{k\>n}\ (k\in\{1,2\})$ is induced by $t_n$, and
 $b_{k\>n}$ will be defined later:
\begin{equation}\label{diagram}
\def\1{u_{n*}u_n^!q_1^!h_1^*f_1^! }
\def\2{u_{n+1*}u_{n+1}^!q_1^!h_1^*f_1^!}
\def\3{u_{n*}q_{1\>n}^!h_{1\>n}^*f_{1\>n}^!y_n^!}
\def\4{u_{n+1*}q_{1\>n+1}^!h_{1\>n+1}^*f_{1\>n+1}^!y_{n+1}^! }
\def\7{u_{n*}q_{2\>\>n}^!h_{2\>\>n}^*f_{2\>\>n}^!y_n^!}
\def\8{ u_{n+1*}q_{2\>\>n+1}^!h_{2\>\>n+1}^*f_{2\>\>n+1}^!y_{n+1}^!}
\def\9{u_{n*}u_n^!q_2^!h_2^*f_2^! }
\def\0{u_{n+1*}u_{n+1}^!q_2^!h_2^*f_2^! }
\begin{CD}
\1 @>c_{1\>n}^{}>>\2 \\
@Vu_{n*}(s_{1\>n})V  V@V V u_{n+1*}(s_{1\>n+1}) V\\
\3 @>b_{1\>n}>>\4 \\
@Vu_{n*}(r_n)V V @V Vu_{n+1*}(r_{n+1})V \\
\7 @>>\under{1.4}{b_{2\>\>n}}>\8 \\
\vspace{-20pt}\\
@Vu_{n*}(s_{2\>\>n}^{-1})V V @V Vu_{n+1*}(s_{2\>\>n+1}^{-1})V \\
\9 @>>\under{1.4}{c_{2\>\>n}^{}}>\0 
\end{CD}
\end{equation}

Let us deal first with the top subrectangle of \eqref{diagram}. (The bottom one is
essentially the same.) To lighten notation, we set $m\!:=n+1$.

Consider the following expansion of the left side of diagram~\eqref{hex},
where all occurrences of ``1" in a subscript have been hidden, and where all the
vertical arrows represent natural closed immersions, so that for each 
$\xi\in \{u,v, x, y\}$, $\xi_n = \xi_m \smcirc \xi_m^n$.
$$
\begin{CD}
U_n @>q_n>> V_n @>h_n>> X_n @>f_n>> Y_n \\
@V u_m^n VV @V v_m^n VV @VV x_m^n V @VV y_m^n V \\
U_m @>q_m>> V_m @>h_m>> X_m @>f_m>> Y_m \\
@V u_m VV @V v_m VV @VV x_m V @VV y_m V \\
\U @>>\under{1.4}q> \V @>>\under{1.4}h> \X @>>\under{1.4}f> \Y 
\end{CD}
$$
The squares in the middle are fiber squares, to which are
associated base-change isomorphisms of the form $h^*x^!\iso
v^!h^*$ (with appropriate subscripts attached). With each $\beta$ indicating
the use of such a base-change isomorphism, and $a_n$ the natural composition
$$
 u_{n*}q_n^!v_{m}^{n\,!}
\iso u_{m*}u_{m*}^{n}u_{m}^{n\>!}q_{m}^{!}
\to u_{m*}q_{m}^!\>,
$$
one sees then that the top rectangle in \eqref{diagram} expands naturally as
\begin{small}
$$
\def\0{u_{n*}u_n^!q^!h^*f^!}
\def\1{u_{m*}u_{m*}^{n}u_{m}^{n\>!}u_{m}^{!}q^!h^*f^!}
\def\2{u_{m*}u_{m}^!q^!h^*f^!}
\def\3{u_{n*}q_n^!v_n^!h^*f^!}
\def\4{u_{n*}q_n^!v_{m}^{n!}v_{m}^!h^*f^!}
\def\5{u_{m*}q_{m}^!v_{m}^!h^*f^!}
\def\6{u_{n*}q_n^!h_n^*x_n^!f^!}
\def\7{u_{n*}q_n^!h_n^*x_{m}^{n!}x_{m}^!f^!}
\def\8{u_{n*}q_n^!v_{m}^{n!}h_{m}^*x_{m}^!f^!}
\def\9{u_{m*}q_{m}^!h_{m}^*x_{m}^!f^!}
\def\ten{u_{n*}q_n^!h_n^*f_n^!y_n^!}
\def\lvn{u_{n*}q_n^!h_n^*x_{m}^{n\>!}f_{m}^!y_{m}^!}
\def\twv{u_{n*}q_n^!v_{m}^{n!}h_{m}^*f_{m}^!y_{m}^!}
\def\thr{u_{m*}q_{m}^!h_{m}^*f_{m}^!y_{m}^!}
\def\vvv{{\text{via } a_n}}
\minCDarrowwidth=.25in
\begin{CD}
\0 @.\Rarrow{12.5em}@.  \1 @>>> \2 \\
@VVV @. @VVV @VVV\\
\3 @.\Rarrow{12.5em}@.  \4 @>\vvv>> \5 \\
@V\beta V\mspace{189 mu}(*)V @. @VV\beta V @VV\beta V\\
\6@>>> \7 @>\beta >> \8 @>\vvv >> \9 \\
@VVV @VVV @VVV @VVV\\
\ten@>>> \lvn @>>\under{1.4}{\beta} > \twv @>\vvv >> \thr 
\end{CD}
$$
\end{small}%
with $b_{1\>n}$ in \eqref{diagram} defined to be the composition 
of the maps in the bottom row. 

It remains only to check commutativity of each of the subrectangles, which is 
a straightforward exercise requiring only the simplest formal properties of
functoriality and pseudofunctoriality,%
\footnote{in particular, ``pseudofunctorial associativity": if
$\gamma_{\psi,\>\varphi}\colon( \varphi\psi)^!\iso\psi^!\varphi^!$ 
is\vspace{.6pt}
the canonical isomorphism then for any composition $\varphi\psi\chi$, 
it holds that
$(\chi^!\gamma_{\psi,\>\varphi})\circ\gamma_{\chi,\>\varphi\psi}$ and
$\gamma_{\chi,\>\psi}(\varphi^!)\circ\gamma_{\psi\chi,\>\varphi}$ are 
\emph{the same}
isomorphism from $(\varphi\psi\chi)^!$ to $\chi^!\psi^!\varphi^!$.} 
except for the subrectangle marked
$(*)$, where one uses the transitivity of flat base change \cite[Lemma
7.5.2(b)]{AJL}.

\smallskip
As for the middle subrectangle in \eqref{diagram}, after noting that
$u_{n*}=u_{m*}u^n_{m*}$ and $y_n^!\cong y_m^{n\>!}y_m^!$ one can ``factor out"
$u_{m*}$ and $y_m^!$, and  then use the rather simple
duality isomorphism for closed immersions, 
under which the natural map $u_{m*}^{n}u_{m}^{n\>!}\G\to\G$ corresponds
to the identity map of $u_{m}^{n\>!}\G$)
$$
{\Hom}_{\D(U_{m})}(u^n_{m*}\E\<,\>\G)\cong
{\Hom}_{\D(U_n)}(\E\<,\>u_m^{n\>!}\G)
\quad \bigl(\E\in\Dqc(U_n),\;\G\in\D(U_{m})\bigr)
$$
to reduce the commutativity question to that for a  diagram of
isomorphisms of functors from
$\Dqcp(Y_m)$ to $\Dqc(U_n)$:
$$
\def\1{q_{1\>n}^!h_{1\>n}^*f_{1\>n}^!y_m^{n\>!}}
\def\2{u_m^{n\>!}q_{1\>m}^!h_{1\>m}^*f_{1\>m}^!}
\def\3{q_{2\>\>n}^!h_{2\>\>n}^*f_{2\>\>n}^!y_m^{n\>!}}
\def\4{ u_m^{n\>!}q_{2\>\>m}^!h_{2\>\>m}^*f_{2\>\>m}^!}
\begin{CD}
\1 @>\Iso>>\2 \\
@V\simeq V V  @V V\simeq V \\
\3 @>\Iso>>\4 
\end{CD}
$$
In this diagram only maps between ordinary schemes appear, so as
before one can identify $h^*$ with $h^!\<$ and define the vertical arrows via
pseudofunctoriality. The horizontal arrows involve flat base change. However, under
the  identification of
$h^*$ with $h^!$ (when $h$ is an open immersion of ordinary schemes), a
base-change isomorphism like $h_{1\>n}^*x_{1\>m}^{n\>!}\iso
v_{1\>m}^{n\>!}h_{1\>m}^*$  becomes identified
with\vspace{1pt} the pseudofunctoriality isomorphism 
$h_{1\>n}^!x_{1\>m}^{n\>!}\iso v_{1\>m}^{n\>!}h_{1\>m}^!$. (This is not trivial,\vspace{1pt}
but is contained in the construction of the pseudofunctor $^{\boldsymbol !}$.)
With this in mind one finds again that pseudofunctoriality
yields the desired commutativity.

This completes the proof of Lemma~\ref{fix}.
\end{proof}

\begin{prop}[\cite{AJL}, 8.3.1]
Let there be given a commutative diagram
\begin{diagram}[size=2em]
    \U            & \rTo^{\st i_1} & \X_1            \\
   \dTo^{\st i_2} &                & \dTo_{\st f_1}  \\
    \X_2          & \rTo_{\st f_2} & \Y
\end{diagram}
of noetherian formal schemes, with\/ $f_1$ and\/~$f_2$
pseudo\kern.6pt-proper and
$i_1$ and\/~$i_2$ open immersions. Then one can define isomorphisms
\[
i_1^*f_1^!\F \iso i_2^*f_2^!\F\>,\qquad\quad
i_1^*\ush{f_1}\F \iso i_2^*\ush{f_2}\F
\qquad\bigl(\F \in \wDqcp(\Y)\bigr)
\]
such that the induced  homology isomorphisms are canonical and functorial.
\end{prop}

\begin{proof}
This is the particular case $q_k={}$identity, $h_k = i_k$, of Lemma~\ref{fix}.
\end{proof}

\begin{prop}[\cite{AJL}, 8.3.2]
   If\/ $f\colon\X\to\Y$ is a pseudo\kern.6pt-proper map of
noetherian formal schemes then
\[
\ush{f}\bigl(\Dc^+(\Y)\bigr)\subset \Dc^+(\X).
\]
\end{prop}

\begin{proof}
As in \emph{loc.\,cit.}\ we may assume that $\Y$ is
affine, say $\Y = \spf(A)$, and that $\X$~can be covered by open subsets
$j : \U \to \X$ such that  $f|_{\U} \!:= j \smcirc \<f\>$ factors as
\[\U \overset{i}{\lto} \spf(B) \overset{h}{\lto} \P \overset{p_{\<1}^{}}{\lto}
\spf(A) \]
where $i$ is a closed immersion, $h$ is an open immersion and $\P$
is the completion of the projective space $\PP^n_{\!\!A}$  along
some closed subset. Now, for $\F \in \Dc^+(\Y)$, Lemma \ref{fix} provides
an isomorphism
$ j^*\<\<\ush{f}\F \cong \ush{i}\>h\<^*\ush{p_1}\F$,
giving a reduction to the two cases \mbox{(a) $f=p_1^{}$}   and (b) $f$ a closed
immersion, cases dealt with at the end of the proof in
\emph{loc.\,cit.}
\end{proof}

\pagebreak[3]

\end{document}